\documentclass[journal,10pt]{IEEEtran}

\usepackage{amsmath,epsfig,cite,amsfonts,amssymb,psfrag}
\usepackage[caption=false]{subfig}
\usepackage{graphicx}
\usepackage{fixltx2e}
\usepackage{pstricks,calc,epsf,xcolor,xspace,dsfont}
\usepackage{color}

\def\no{\nonumber}

\newtheorem{theorem}{Theorem}

\newtheorem{rem}{Remark}
\newtheorem{defi}{Definition}

\begin{document}

\title{Sequential Detection of Deception Attacks in Networked Control Systems with Watermarking}
\vspace{-1.7cm}
\author{Somayeh Salimi, Subhrakanti Dey and Anders Ahl\' en \thanks{The authors are with the Division of Signals and Systems, Uppsala University, Box 534, Uppsala, SE-75121, Sweden (e-mail: subhrakanti.dey@signal.uu.se)}%
}

\maketitle
\thispagestyle{empty}
\pagestyle{empty}

\maketitle
\vspace{-1.7cm}

\begin{abstract}
In this paper, we investigate the role of a physical watermarking signal in quickest detection of a deception attack in a scalar linear control system where the sensor measurements can be replaced by an arbitrary stationary signal generated by an attacker. By adding a random watermarking signal to the control action, the controller designs a sequential test based on a Cumulative Sum (CUSUM) method  that accumulates the log-likelihood ratio of the joint distribution of the residue and the watermarking signal (under attack) and the joint distribution of the innovations and the watermarking signal under no attack. As the average detection delay in such tests is asymptotically (as the false alarm rate goes to zero) upper bounded by a quantity inversely proportional to the Kullback-Leibler divergence(KLD)  measure between the two joint distributions mentioned above, we analyze the effect of the watermarking signal variance on the above KLD. We also analyze the increase in the LQG control cost due to the watermarking signal, and show that there is a tradeoff between quick detection of attacks and the penalty in the control cost. It is shown that by considering a sequential detection test based on the joint distributions of residue/innovations and the watermarking signal, as opposed to the distributions of the residue/innovations only, we can achieve a higher KLD, thus resulting in a reduced average detection delay. Numerical results are provided to  support our claims.
\end{abstract}
%
%
\vspace{-.1cm}

\section{Introduction}
Attacks on cyber-physical systems (CPS) can affect the integrity, availability and confidentiality in CPS.  Examples range from deception based attacks such as  {\em false-data-injection} \cite{kallesec15}, sensor and actuator attacks, {\em replay attacks}, and also {\em denial of service attacks} \cite{mo_sinopoli_15} on the underlying networked control system (NCS).
Deception attacks refer to scenarios where integrity of control packets or measurements are compromised by altering the behaviour of sensors and actuators. In particular, false data injection attacks are introduced by injecting incorrect or misleading measurements or control inputs. Replay attacks are carried out by hijacking the sensors, recording the sensor measurements for a period of time, and then repeating such measurements to the controller while injecting a harmful exogenous signal into the system. On the other hand, denial of service attacks can be carried out by an adversary compromising the availability of resources to the CPS, e.g., by jamming the communication channel. 
Documented defence mechanisms can range from attack identification and detection, intrusion detection as well as physical watermarking of valid control signals. Most of these defence mechanisms have been developed to tackle specific types of attacks, whereas a generalized unified approach for attack identification and detection 
is developed by adopting a descriptor system modelling framework for CPS \cite{bullo_sec_15} and applications illustrated for power and water networks. Linear state estimation with corrupted measurements has been also studied in 
\cite{tabuada_sec_14} where the maximum number of faulty sensors is characterized and a decoding algorithm for detecting corrupted measurements is presented. 

The defence mechanism of relevance to this paper is the idea of {\em physical watermarking of control signals}. Traditionally, digital watermarking has been used extensively in audio and image processing for authentication purposes, where a specific signal is embedded in the transmitted message/document, and is later used to identify the rightful owner of the message. The idea of  physical watermarking in NCS is similar, where a random signal is added to the control signal, and under normal operations, the effect of this watermarking signal should be present in the system output. However, when the system is attacked or compromised and sensor measurements are substituted by injection of false data, the expected effect of the watermarking signal will be absent or perturbed, thus leading to  a statistical test which can detect the presence of an attacker. Two most recent works that deal with design and analysis of physical watermarking for NCS are \cite{mo_sinopoli_15} and \cite{kumar}. In \cite{mo_sinopoli_15}, the authors consider a linear state space model under a replay attack, and design an optimal watermarking signal (added to the true control signal) which maximizes the Kullback-Leibler Divergence (KLD) measure between the densities of the residual before and after the attack, subject to a constraint on the loss of linear quadratic (LQ) control cost due to the addition of the watermarking signal. A typical 
$\chi^2$ - failure detector is used to detect an attack when a watermarking signal is added. In \cite{kumar}, 
the authors proposed a model where the attacker also replaces the sensor measurements by its own simulated signal, which tries to mimic the nominal system  without the knowledge of the watermarking signal in the control input. The key result in \cite{kumar} develops two tests at the actuator that the attacker has to pass to remain stealthy, but this is only possible if the attacker replaces the true sensor outputs by a signal of zero average energy. 

In detecting attacks in CPS, it is of paramount importance that attack detection happens with minimum delay, thus favouring {\em quickest sequential detection} based methods. The watermark design techniques employed in 
\cite{mo_sinopoli_15, kumar} are not designed specifically for this purpose, and the statistical detection tests developed in \cite{kumar} are {\em asymptotic} in nature, thus relying on collecting a large number of system outputs in practice. In this paper, we will therefore focus on design and analysis of physical watermarking signals that minimize the average detection delay in sequential detection methods, while still keeping the system performance within a prescribed safety limit - as demanded by resilience requirements of CPS under attacks
\cite{cardenas_sastry_08}. 

In particular, we consider a scalar networked linear control system, where the attacker launches a deception attack at a certain unknown but deterministic time point, by injecting a false measurement sequence that replaces the true sensor measurements. The estimator/controller employs standard optimal linear quadratic Gaussian (LQG) control based on the received measurement sequence without knowing whether there has been an attack or not. In order to aid the detection of the attacker, which on the other hand tries to remain stealthy, the controller adds a random watermarking signal to the control signal, which is only known to the controller/actuator, and not the adversary.  Furthermore, the controller employs a sequential detection test based on the cumulative sum (CUSUM) algorithm that is well known to minimize the average detection delay under a constraint on the mean time between false alarms. This sequential test is based on the log-likelihood ratio of the joint distribution of the residue (measurement prediction error) and the watermarking signal before and after the attack. Since an asymptotic (as the false alarm rate goes to zero) upper bound on the average detection delay  is inversely proportional to the Kullback-Leibler divergence (KLD) measure between these joint  distributions before and after the attack, we analyze the behaviour of the KLD measure with respect to the variance of the watermarking signal. While increasing the watermarking signal variance increases the KLD, it also increases the corresponding LQG control cost. The behaviour of the increase in the control cost due to the watermarking signal is also analyzed, illustrating the tradeoff between quickest detection, and the penalty in the control cost. Unlike previous works which consider KLD between the distributions of the residue signal (under attack) and the innovations (before the attack) only, we show that by considering the joint distributions of the residue/innovation and the watermarking signal, we can increase the KLD even further, thus reducing the average detection delay. Numerical results confirm our  findings. 

\section{Problem Setup}
\subsection{System model}
\noindent We consider the following 
architecture of a networked control system.
 The single-input single-output, linear time invariant system is
modeled as:
\begin{equation}
x_{k+1}=Ax_{k}+Bu_{k}+w_{k}\label{syseq}
\end{equation}
in which $x_{k}\in \mathbb{R}$  is the state variable and $u_{k}\in
\mathbb{R}$ is the control input at time $k$ generated by the
controller. $w_{k}\in \mathbb{R}\sim\mathcal{N}(0,Q)$ is the process
noise at time $k$ which is assumed to be an independent and identically distributed ({\em i.i.d}) random process.
A sensor reports its ({\em scalar}) observations to the controller  in the
following form at time $k$:
\begin{equation}
y_{k}=Cx_{k}+v_{k}\label{ydef}
\end{equation}
in which $v_{k}\sim \mathcal{N}(0,R)$ is i.i.d. measurement noise
that is independent of the process noise $w_{k}$. Note that although we consider a scalar state-space system, we still use uppercase letters for the system parameters $A, B, C, Q, R$, with a slight abuse of notation. We assume that the
system has started at time $t = -\infty$ and currently is in
steady-state condition, as stabilizability and detectability are guaranteed for a scalar system. Then the optimal state estimate equations, based
on Kalman filtering, are given as
\begin{eqnarray}
\hat{x}_{k+1|k}=A\hat{x}_{k|k}+Bu_{k}\\
\hat{x}_{k|k}=\hat{x}_{k|k-1}+K\gamma_{k}\label{xestimatedef}
\end{eqnarray}
where $\hat x_{k+1|k} = E[x_{k+1} |{\cal Y}_k]$, and $\hat x_{k|k} = E[x_k | {\cal Y}_k]$ are the Kalman predicted and filtered state estimate, respectively based on received measurements up to time $k$, given by ${\cal Y}_k$. Also, 
\begin{equation}
K=\frac{CP}{C^{2}P+R}
\end{equation}
is the steady-state Kalman gain
and where $P$ is the steady-state minimum mean-squared error (MMSE) estimation error variance
$E(x_{k}-\hat{x}_{k|k-1})^{2}$ obtained from the solution to the algebraic Riccati
equation
\begin{equation}
P=A^{2}P+Q-A^{2}C^{2}P^{2}(C^{2}P+R)^{-1}.\label{MMSE}
\end{equation}

In \eqref{xestimatedef} the innovation sequence $\gamma_{k}$ is the  defined as

\begin{equation}
\gamma_{k}\triangleq y_{k}-C\hat{x}_{k|k-1}.\label{gammadef}
\end{equation}
We assume that the sensor is connected to the estimator/controller via a link that is susceptible to malicious attacks.
In the system equation \eqref{syseq}, the control signal $u_{k}$ is sent
by the controller (which  is assumed to be co-located with the actuator), to the sensor as a
linear function of the filtered state estimate, such that
$u_{k}=f(\hat{x}_{-\infty}^{k})$ minimizes the infinite-horizon
LQG cost:
\begin{equation}
J=\lim_{T\rightarrow\infty}E\frac{1}{2T+1}[\sum\limits_{k=-T}^{T}(Wx_{k}^{2}+Uu_{k}^{2})]\label{LQG}
\end{equation}
where $W$ and $U$ are positive weights. The LQG control policy
results in a fixed-gain linear control signal as

\begin{equation}
u_{k}=L\hat{x}_{k|k}\label{contsig}
\end{equation}
in which $L$ is given by
\begin{equation}
L=\frac{-ABS}{B^{2}S+U}
\end{equation}
and $S$ is the solution obtained from the algebraic Riccati equation
\begin{equation}
S=A^{2}S+W-A^{2}B^{2}S^{2}(B^{2}S+U)^{-1}.
\end{equation}
\subsection{Attack model}
We assume that the adversary can launch an attack against the integrity
of the sensor measurements such that the estimator/controller, instead
of receiving the true measurement, $y_{k}$ sent by the honest sensor,
receives $z_{k}$, which is injected by the attacker. 
Furthermore, we assume that
the attacker knows the system parameters $A, B,
C, Q$ and $R$ and also the control policy, i.e., $L$ but not
necessarily the true sensor measurements $y_{k}$.
On the other hand, we assume that the control
signal is not tampered with by the adversary.

The goal of the attacker is to change the performance of the control
system by sending fake observations, $z_{k}$, that replaces the true ones and while doing so
remain undetected. In the
general form of the attack, we assume that the attack is injected into the system at
time $k$, i.e., instead of $[y_{k}, y_{k+1}, y_{k+2},... ]$,
$[z_{k}, z_{k+1}, z_{k+2},... ]$ is received by the controller.
It is easy to show that if the attacker only modifies the sensor measurement at time $k$, and left the subsequent measurements undisturbed, or modifies the sensor measurements only after sufficiently long  intervals, the attacker can remain undetected as the difference between the true innovations and the false innovations based on the altered measurement sequence will go to zero exponentially fast. On the other hand, the effect of the attack at one single point in time will also be forgotten exponentially due to stabilizability and detectability properties of the control system which are automatic for the scalar case when the optimal Kalman filter and LQG controller are applied. This is obviously not useful from the attacker's point of view. Hence, in the following section, we consider an attack model where the attacker continuously replaces the true measurement $y_m$ by a fake measurement $z_m$ for all $m \geq k$. {\black This is a kind of spoofing attack which can be accomplished by the adversary, even without having access to $y_k$ itself, by jamming or overpowering the true sensor signal if sent over wireless. However, if the sensor signal is not sent over wireless the adversary might be able to hijack it in another way and replacing the $y_k$ with $z_k$ in a so called man-in-the-middle attack. Most protocols used today would not be able to detect such an attack. Nevertheless, the objective of the attacker is to remain stealthy for a sufficiently long period of time over which the attack takes place, to cause maximum damage to the control system.}

In this paper, we will assume that the attacker does not need to know the true sensor measurements but can simply alter them by injecting (as we will assume for the rest of this paper) the sequence
$\{z_{k}\}$, which is stationary with statistics
\begin{eqnarray}
E(z_{k}^{2})=\sigma_{z}^{2},\\
E(z_{k},z_{k-k'})=\rho^{k'}\sigma_{z}^{2}\label{zcorrelation}
\end{eqnarray}
in which $\rho<1$. Depending on whether the attacker physically compromises the sensor node or simply replaces the sensor measurements by injecting a strong interfering signal, it may also need to know the encryption algorithm used by the networked control system. However, it is  common to assume that the adversary has full knowledge of all system parameters and protocols, as is often done in cryptography according to the notion of ``security through obscurity'' known as Kerckhoffs's principle, or also according to {\em Shannon's maxim}, which essentially assumes that 
``the enemy knows the system.'' {\black Knowing the system is a sensible assumption since then the adversary can cause maximum damage, a situation that is essential to detect as fast as possible. }

\subsection{Attack stealthiness}
To determine whether an attack is present in the control system or not we shall rely on a hypothesis
testing procedure based on  the following two hypotheses: \vspace{.25cm}

$H_{0}$: No attack (the controller receives the true sequence $y_{k}$)
\vspace{.5mm}

    $H_{1}$: Attack (the controller receives the false sequence $z_{k}$)
\vspace{.25cm}

Let $p_{k}^{F}$ represent the false alarm probability, i.e., deciding $H_{1}$
when $H_{0}$ is true and let $p_{k}^{D}$ represent the detection probability, i.e., deciding $H_{1}$ when $H_{1}$ is true, at time $k$. Furthermore, define $\widetilde{\gamma}_k$ to be the innovation signal $z_k - C \hat x^F_{k|k-1}$, where 
$ \hat x^F_{k|k-1}$ is the inaccurate Kalman predictor designed in the presence of an attack based on the received sequence $\{z_k\}$. 
Let $\widetilde{\gamma}_{1}^{k}$ and $\gamma_{1}^{k}$ represent the sequences $\{\widetilde{\gamma}_j\}_{j=1}^{k}$ $\{\gamma_j\}_{j=1}^{k}$, respectively.
The goal is to design a detector which, with high probability can detect an attack while keeping the false alarm probability as small as possible. It is common to design a hypothesis testing procedure that decides in favour of $H_0$ or $H_1$ 
based on testing the innovation sequence $\widetilde{\gamma}_{1}^{k}$ (under attack) and the true innovation sequence 
$\gamma_{1}^{k}$. 

In detection theory, the performance of the detector can be characterized by the trade-off between $p_{k}^{F}$ and $p_{k}^{D}$. Following \cite{stealthyattack},\cite{enoch}, we introduce the following definition of a stealthy attack:

\vspace{.2cm}
\begin{defi}
For $\epsilon>0$ and $0<\delta<1$, an attack is $\epsilon$-stealthy
if for any detector that satisfies $0<1-p_{k}^{D}\leq\delta$, it holds that
\begin{equation}
\limsup_{k\rightarrow\infty}-\frac{1}{k}\log(p_{k}^{F})\leq\epsilon\label{stealthy}
\end{equation}
\end{defi}

It was shown in \cite{enoch} that condition \eqref{stealthy} is
equivalent to
\begin{equation}
\limsup_{k\rightarrow\infty} \frac{1}{k}D(f_{\widetilde{\gamma}_{k}}\|f_{\gamma_{k}})\leq\epsilon
\end{equation}
when the hypothesis $H_0$ for no attack assumes the innovation sequence $\gamma_{1}^{k}$, and the residues
$\widetilde{\gamma}_{1}^{k}$ for $H_1$.   Here, $D(f_{\widetilde{\gamma}_{k}}\| f_{\gamma_{k}})$
is the Kullback-Leibler Divergence (KLD) between the sequences
$\widetilde{\gamma}_{1}^{k}$ and $\gamma_{1}^{k}$ defined as:

\begin{equation}
D(f_{\widetilde{\gamma}} \| f_{\gamma})=\int_{-\infty}^{\infty}f_{\widetilde{\gamma}}(\gamma_{1}^{k})\log\frac{f_{\widetilde{\gamma}}(\gamma_{1}^{k})}{f_{\gamma}(\gamma_{1}^{k})}d\gamma_{1}^{k}\ .
\end{equation}
where $f_{\widetilde{\gamma}}, f_{\gamma}$ are the (stationary) distributions of the sequences 
$\{\widetilde{\gamma}_{k}\}$ and $\{\gamma_{k}\}$, respectively.  Clearly, for a given $\epsilon$ designed by the attacker, 
the objective for the control system designer is to detect the attacker, and hence increase the value of the 
quantity $D = \limsup_{k\rightarrow\infty} \frac{1}{k}D(f_{\widetilde{\gamma}_{k}}\|f_{\gamma_{k}})$, an expression for which was provided in \cite{enoch}. This leads us to the next section, where we employ a physical watermarking mechanism to increase an appropriate KLD measure based on the joint distributions of the innovations/residues and the random watermarking signal, thus making it difficult for the attacker to remain undetected through a sequential detection test designed accordingly. Note that we do not discuss how the adversary designs $\epsilon$ in response to the sequential detection test employed by the control system designer in this paper, a topic which will be further investigated in a game theoretic setting in future work.

\subsection{Defence mechanism based on physical watermarking}
As explained above, the attacker can choose an intelligent policy to inject  false observations and tries to remain undetected. This however relies on the the fact that the control system is influenced by process  and measurement noises, which produce uncertainty in favour of the
attacker. 

To protect the system against these active attacks, a key idea is to add a random watermarking  signal, known only by the controller (and not to the attacker, although the attacker may know the statistics of the watermarking signal), to the
control sequence $u_{k}$. In particular, the controller adds the watermarking sequence $e_{k}$ to the control signal,
i.e., 
\begin{equation}
u_{k}=L\hat{x}_{k}+e_{k}\label{uwithw}
\end{equation}
where $e_{k}$ is assumed to be an i.i.d. zero-mean Gaussian sequence with variance
$\sigma_{e}^{2}$. The idea of adding such a physical watermarking signal was proposed in \cite{mo_sinopoli_15} in the context of detecting {\em replay attacks}, and further extended and analyzed in the context of dynamic watermarking 
in \cite{kumar}. In general, the signal $e_{k}$ can be a stationary Gauss-Markov process as shown in 
\cite{mo_sinopoli_15}, although for the purpose of this paper, we assume it to be i.i.d.

By adding this sequence the controller is provided  with a tool to check if the received signal from the sensors bear any correlation with the watermarking sequence or not. If the attacker injects a false observation $z_{k}$, which is naturally independent of the watermarking signal, then this can be detected by the controller, even though the attacker may know the statistics of the
watermarking signal.

In \cite{mo_sinopoli_15}, a $\chi^2$ failure detector based on the residue vector (which is either $\gamma_k$ or 
${\widetilde{\gamma}_k}$) was suggested for detecting an attack, whereas in \cite{kumar}, two asymptotic tests were proposed to detect an attack. Both of these schemes require a sufficiently large number of samples to be used for the test in practice to achieve a good detection probability with a constraint on the false alarm rate. It is of course, of utmost importance to detect an attack as soon as possible, and this motivates us to consider a non-Bayesian sequential detection method under the assumption that the attack takes place at a fixed but unknown point of time. In particular, the cumulative-sum (CUSUM) method which minimizes the average detection delay subject to a constraint 
on the mean time between false alarms, also known as Lorden's method \cite{quickest_detection_book}. However, instead of comparing directly the distribution of $\gamma_k$ and
${\widetilde{\gamma}_k}$, we propose a  detection mechanism as follows.
The controller, upon receiving the observation $y_{k}$ (which is not known to be the true $y_{k}$ or the false $z_{k}$)
calculates $\gamma_{k}$ (or $\widetilde{\gamma}_k$ )and
computes
\begin{equation}
S_{k}=\max(0,
S_{k-1}+\log\frac{f_{\widetilde{\gamma}_{k}, e_{k-1}}(\widetilde{\gamma}_{k}, e_{k-1})}{\black f_{\gamma_{k}, e_{k-1}}(\widetilde{\gamma}_{k}, e_{k-1})})\ .\label{seqtest}
\end{equation}
where $f_{\widetilde{\gamma}_{k}, e_{k-1}}$ and $f_{\gamma_{k}, e_{k-1}}$ denote the joint distribution between the 
residue signal and the watermarking signal. 
The controller then decides on ``attack" or ``no attack" based on
the following policy:
\begin{eqnarray}
\text{The system is under attack if} \ \ \ S_{k}>\alpha,\no\\
\text{The system is not under attack if} \ \ \ S_{k}<\alpha\label{detect1}
\end{eqnarray}
where $\alpha\triangleq |\log p^{F}|$ .

The above policy can be justified in the way that if the received
observation by the controller is the true one, then
\begin{align}
\gamma_{k} & =y_{k}-C\hat{x}_{k|k-1}\nonumber \\
&=Cx_{k}+v_{k}-C(A+BL)\hat{x}_{k-1|k-1}-CBe_{k-1}\no\\
&=CA(x_{k-1}-\hat{x}_{k-1|k-1})+Cw_{k-1}+v_{k}\ .\label{eq1}
\end{align}

meaning that $\gamma_{k}$ is uncorrelated
with the watermarking signal $e_{k-1}$.  On the
contrary, if the received observation by the controller is the false
$z_{k}$, then
\begin{align}
 \tilde{\gamma}_{k} & =z_{k}-C\hat{x}^F_{k|k-1} \nonumber \\ 
& =z_{k}-C(A+BL)\hat{x}^F_{k-1|k-1}-CBe_{k-1}.\label{eq2}
\end{align}
Thus, it is evident that the false innovations $\tilde{\gamma}_{k}$ is correlated with watermarking
signal $e_{k-1}$ and we can conclude that the control system is under attack.

It is worth mentioning that  one might be tempted to conduct a sequential test 
 based on the log-likelihood ratio 
$\log\frac{f_{\widetilde{\gamma}_{k}}(\widetilde{\gamma_{k}})}{f_{\gamma_{k}}(\gamma_{k})}$ (i.e, based on the log-likelihood ratio of the distributions of the residue under attack and innovations (no attack)), 
instead of $S_k$ defined in \eqref{seqtest}, which is based on the joint distributions of the residue/innovations and the watermarking signal. In the following sections, we will illustrate how our suggested test quantity $S_k$ can reduce  the average detection delay as opposed to using the log-likelihood ratio based on the residue/innovations only, as mentioned above.

\section{Main results}

To analyze our suggested detection approach further we will use
the Average Detection Delay (ADD) as a measure to quantify performance. It is well known that 
\cite{quick, quickest_detection_book} when the observations before and after the change are i.i.d, it is shown that, as the mean time between false alarms goes to infinity (or false alarm rate $p^{F}$ goes to zero) 
the ADD is asymptotically upper bounded by $\frac{|\log p^{F}|}{I_1}$ where $I_1$ corresponds to the KLD between the distributions after and before the change. Although originally derived for i.i.d. sequences, these asymptotic upper bound results have been extended to the case of dependent but stationary sequences in \cite{tartakovsky_dependent}, which allows to write the following asymptotic upper bound on the ADD for the proposed sequential test based on \eqref{seqtest}, \eqref{detect1}:
\begin{align}
& \frac{|\log p^{F}|}{D(f_{\widetilde{\gamma}_{k}, e_{k-1}}\|f_{\gamma_{k}, e_{k-1}})}\ .
\end{align}

Clearly, for  a fixed $p^{F}$, the upper bound on the ADD is inversely proportional to the KLD between the joint distributions before and after the attack. In the following theorem, we obtain an expression for 
$ D(f_{\widetilde{\gamma}_{k}, e_{k-1}}\|f_{\gamma_{k}, e_{k-1}})$
corresponding to our proposed detection approach.

\begin{theorem}
{\black Consider the joint distributions between the watermarking signal and the true and false innovations, respectively, i.e.,  $f_{\gamma_{k}, e_{k-1}}$ and $f_{\widetilde{\gamma}_{k}, e_{k-1}}$.
The KLD between these joint distributions, $ D(f_{\widetilde{\gamma}_{k}, e_{k-1}}\|f_{\gamma_{k}, e_{k-1}})$,  is then given by:}

\begin{align}
D(f_{\widetilde{\gamma}_{k}, e_{k-1}}\|f_{\gamma_{k}, e_{k-1}}) & = \frac{1}{2}\log(\frac{1}{1-\lambda^{2}})
\nonumber \\
& +
\frac{1}{2}(\frac{\sigma_{\tilde{\gamma}}^{2}}{\sigma_{\gamma}^{2}}-1-\log\frac{\sigma_{\tilde{\gamma}}^{2}}{\sigma_{\gamma}^{2}})\label{KLD}
\end{align}

in which
\begin{align}
&\lambda=\frac{-BC\sigma_{e}}{\sigma_{\tilde{\gamma}}}\\
&\sigma_{\tilde{\gamma}}^{2}=[(1-\frac{\rho
CK(A+BL)}{1-\rho\mathcal{A}})^{2} \nonumber \\
& +\frac{(1-\rho^{2})C^{2}K^{2}(A+BL)^{2}}{(1-\mathcal{A}^{2})(1-\rho\mathcal{A})^{2}}]\sigma_{z}^{2}
+\frac{B^{2}C^{2}}{1-\mathcal{A}^{2}}\sigma_{e}^{2}\\
&\sigma_{\gamma}^{2}=C^{2}P+R\label{sigmagammanoattack}
\end{align}
where $ \mathcal{A}\triangleq (1-CK)(A+BL) $ (note also that $\mathcal{A} < 1$ from stabilizability and detectability which is automatic for the scalar case) and $P$ is calculated
according to \eqref{MMSE}. Finally, $|\lambda| < 1$, as shown in the proof.
\end{theorem}

\begin{IEEEproof}
See Section V.A.
\end{IEEEproof}

\vspace{.4cm}
 In Theorem 1, $D(f_{\widetilde{\gamma}_{k}, e_{k-1}}\|f_{\gamma_{k}, e_{k-1}})$ is the difference between
the joint distributions between the innovations and the watermarking signal, i.e., 
$f_{\gamma_{k}, e_{k-1}}(\gamma_{k},e_{k-1})$ and $f_{\tilde{\gamma}_{k}, e_{k-1}}(\tilde{\gamma}_{k},e_{k-1})$, for the healthy and
attacked systems, respectively. With a fixed false alarm probability, we need to
make these distributions as distinguishable as possible to avoid 
unnecessary detection delays. Observing the expression \eqref{KLD}
we note that $D(f_{\widetilde{\gamma}_{k}, e_{k-1}}\|f_{\gamma_{k}, e_{k-1}})$ is a monotonically increasing function of
$\sigma_{e}^{2}$. Hence, an increase in the watermarking signal
power results in a decrease in the ADD due to an increasing KLD. 
However, the decreased ADD comes at a cost: by increasing the watermarking signal
power we also diverge from the optimal LQG cost given by \eqref{LQG}, \eqref{contsig}.
Hence,  there is a tradeoff between reducing the average detection delay 
and the system performance in terms of the increase of the LQG
cost. Note that it can be shown that the difference in the LQG cost due to use of the watermarking signal is higher for the healthy system than that of the system under attack. Hence we elaborate on this issue in the following theorem, by considering the difference in the LQG cost for the healthy system (the worse of the two scenarios). See \cite{mo_sinopoli_15} for a similar treatment.
\begin{theorem}
Consider the LQG cost \eqref{LQG} with weighting factors $W$ and $U$ and let 
$\Delta LQG$ represent the acceptable increase in the LQG cost  from the optimal LQG
cost for the system under no attack. Then, the watermarking signal variance is related to the increase in LQG cost  as follows:
\begin{align}
& \sigma_{e}^{2}=\frac{\Delta
LQG}{U+\frac{B^{2}(W+L^{2}U)}{1-(A+BL)^{2}}}
\end{align}
\end{theorem}

\begin{IEEEproof}
See Section V.B.
\end{IEEEproof}

\vspace{.4cm}
\begin{rem}
As stated in the previous section, instead of the distribution of
the innovations in our detection policy, we considered the joint
distribution between the innovations and the watermarking signal. The
benefit of using the joint distribution instead of using the innovations only can be immediately observed in the expression
of the $D(f_{\widetilde{\gamma}_{k}, e_{k-1}}\|f_{\gamma_{k}, e_{k-1}})$ in \eqref{KLD}. By considering
the innovations only, the corresponding $D(f_{\widetilde{\gamma}_k} \|f_{\gamma_{k}})$ expression equals
the second term in \eqref{KLD} i.e.,
$\frac{1}{2}(\frac{\sigma_{\tilde{\gamma}}^{2}}{\sigma_{\gamma}^{2}}-1-\log\frac{\sigma_{\tilde{\gamma}}^{2}}{\sigma_{\gamma}^{2}})$.
By using the joint distributions we also obtain the first term in  \eqref{KLD}, which is positive, leading to a larger KLD, thus making it more difficult for the attacker to remain stealthy. \end{rem}


\section{Numerical results}
\begin{figure}[htbp]
\centering
\begin{minipage}[t]{\columnwidth}
\includegraphics[height=0.7\columnwidth,width=\columnwidth]{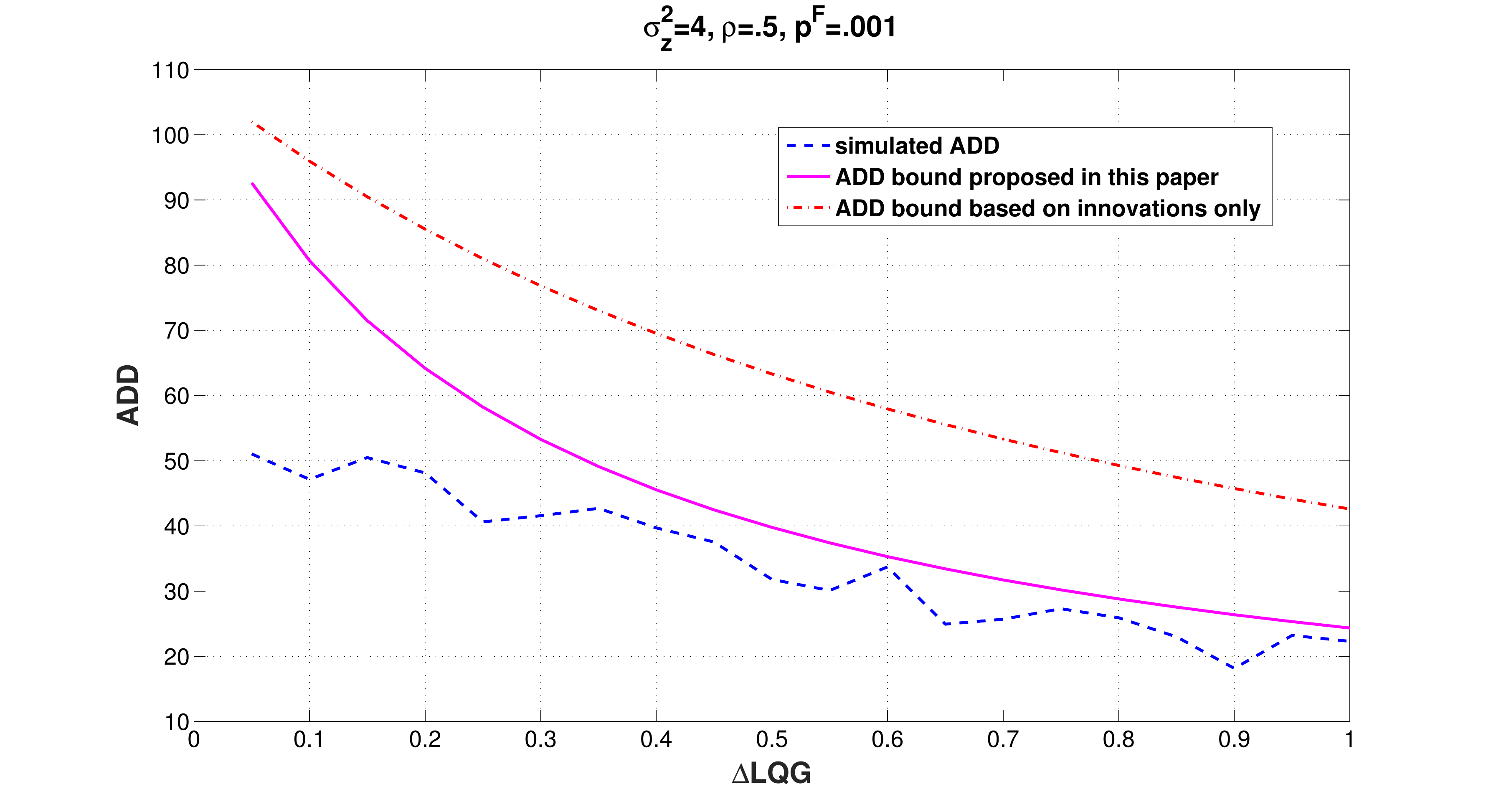} 
\end{minipage}%
\hfill%
\begin{minipage}[t]{\columnwidth}
\includegraphics[height=0.7\columnwidth,width=\columnwidth]{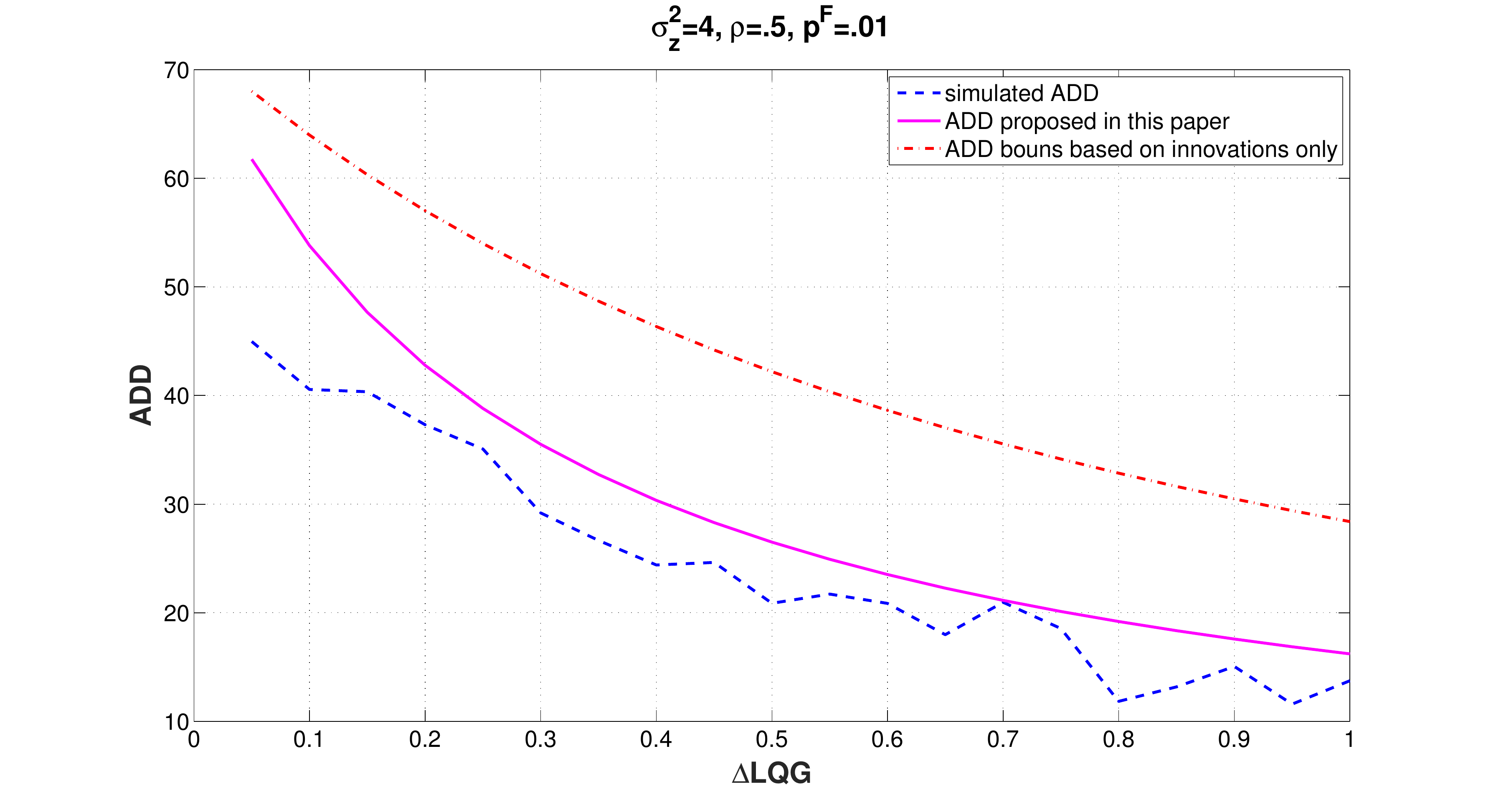}
\end{minipage}

 \caption{\footnotesize{ADD in terms of $\Delta LQG$ for $p^{F}=.001$ and $p^{F}=.01$}}
\label{fig1}
\end{figure}

In this section we will investigate the tradeoff between ADD and $\Delta LQG$.
This is shown in Fig.~\ref{fig1} for two different values of $p^{F}=.001$ and
$p^{F}=.01$. The blue dashed line represents the simulated ADD according to
our detection policy in \eqref{detect1} while the solid pink line
 shows the ADD upper bound
according to our result in Theorem 1.  The red dash-dot line shows the ADD upper bound for  the detection policy 
based on the KLD between the distributions of the residue signal (under attack) and the true innovations (no attack). This is referred  to in the graphs as ``ADD bound based on innovations only''.   It should be noted that in previous works such as \cite{kumar}, some numerical results based on a sequential detection are presented, although the actual detection policy is not clearly stated.   The results in Fig.~\ref{fig1} are depicted for the values $A=.7,
B=C=R=Q=W=1, U=.4$, $\sigma_{z}^{2}=4$, $\rho=.5$ and the real ADD
blue dashed line is calculated based on the average over 100 random realizations of the sequential detection algorithm.
 Comparison between the solid pink line and the red dash-dot
line in Fig. 1  demonstrates the reduction in the ADD upper bounds  due to using our proposed sequential detection test
 compared to a sequential detection based on innovations/residues only.

\begin{figure}[htbp]
\centering
\begin{minipage}[t]{\columnwidth}
\includegraphics[height=0.7\columnwidth,width=\columnwidth]{new3-eps-converted-to.pdf} 
\end{minipage}%
\hfill%
\begin{minipage}[t]{\columnwidth}
\includegraphics[height=0.7\columnwidth,width=\columnwidth]{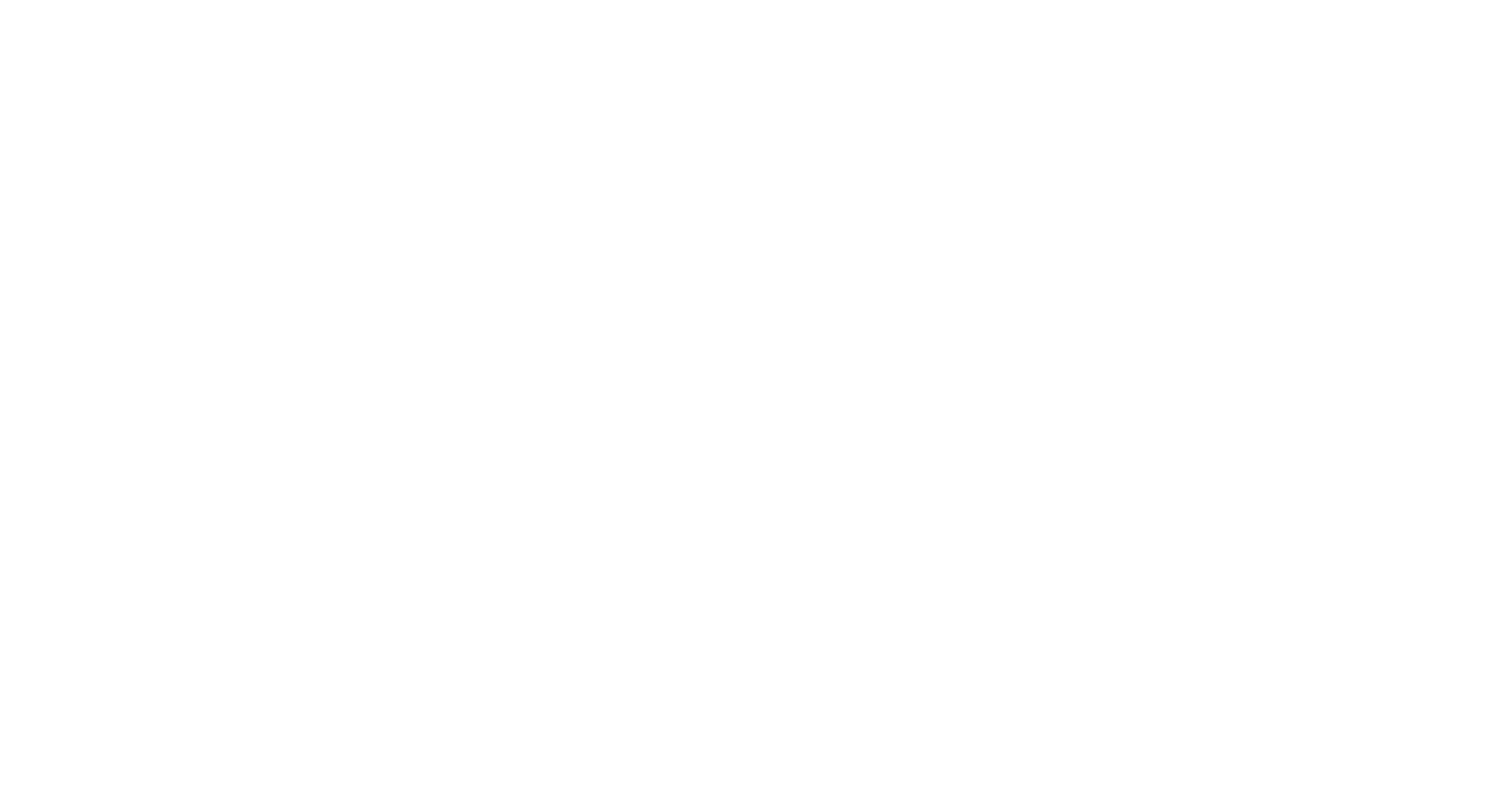}
\end{minipage}
\\
\begin{minipage}[t]{\columnwidth}
\includegraphics[height=0.7\columnwidth,width=\columnwidth]{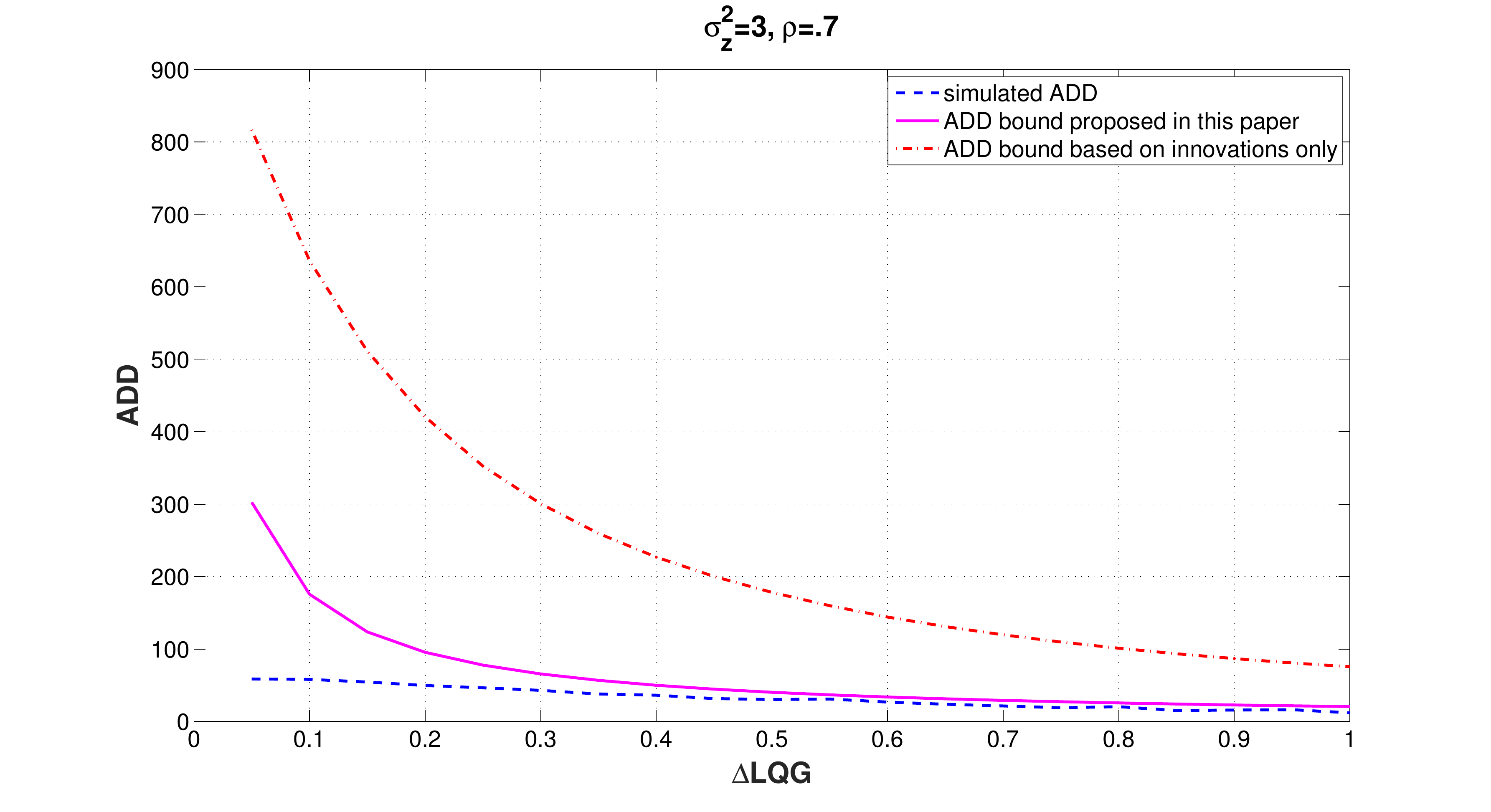} 
\end{minipage}
\hfill%
\begin{minipage}[t]{\columnwidth}
\includegraphics[height=0.7\columnwidth,width=\columnwidth]{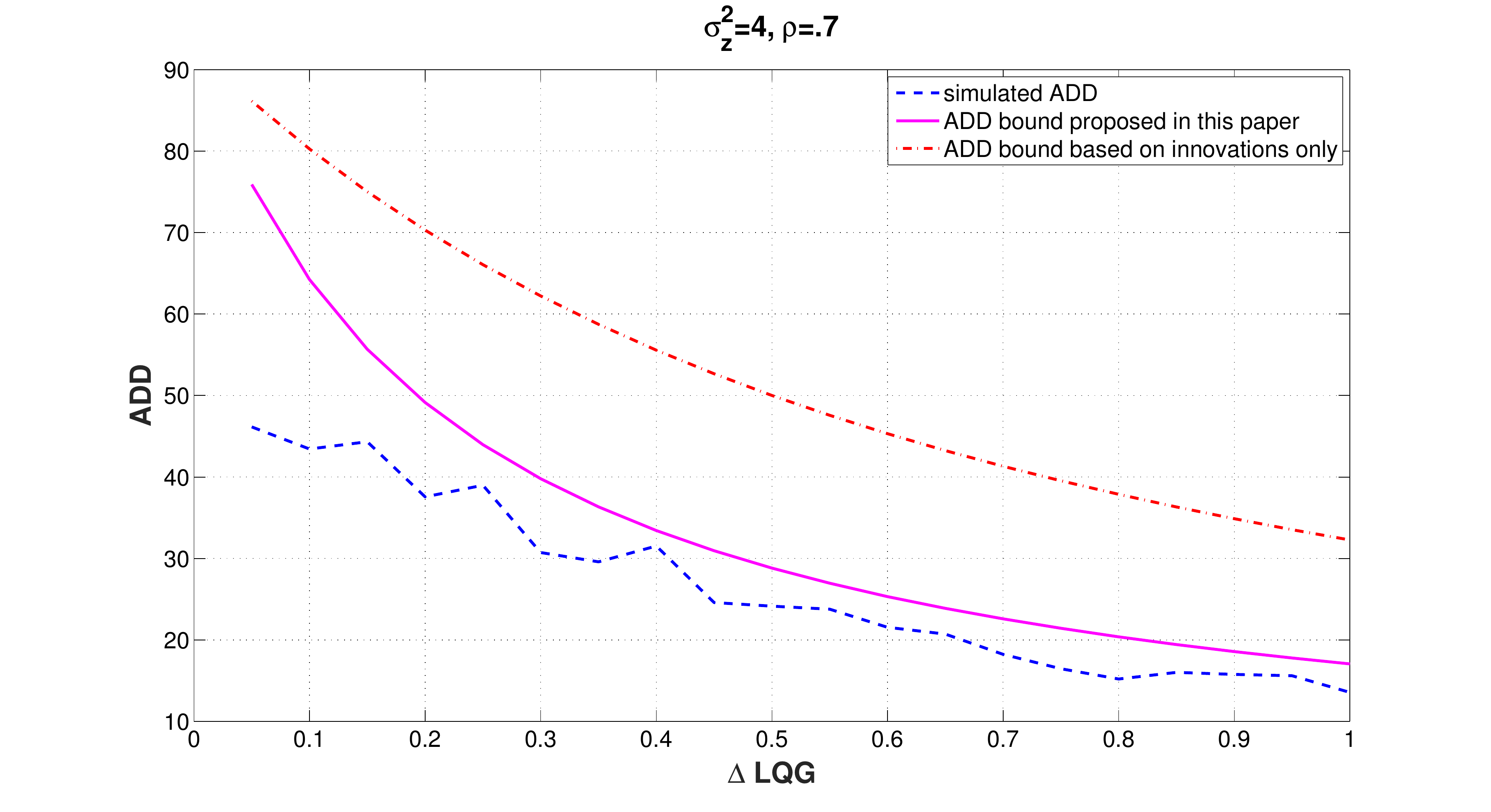}
\end{minipage}
\caption{\footnotesize{\footnotesize{ADD in terms of $\Delta LQG$ for different values of $\sigma_{z}^{2}$ and $\rho$}} }
\label{fig2}
\end{figure}

In Fig.~\ref{fig2} we investigate the effect of different values of
$\sigma_{z}^{2}$ and $\rho$ for fixed value of the parameters $A=.7,
B=C=R=Q=W=1, U=.4,p^{F}=.01$. Intuitively, ADD is a decreasing
function of $\sigma_{z}^{2}$ and an increasing function of $\rho$.
This conforms with the theoretical analysis of the $D(f_{\widetilde{\gamma}_{k}}\|f_{\gamma_{k}})$ in
\eqref{KLD} in which $\frac{\partial D(f_{\widetilde{\gamma}_{k}, e_{k-1}}\|f_{\gamma_{k}, e_{k-1}})}{\partial
\sigma_{z}^{2}}>0$ and $\frac{\partial D(f_{\widetilde{\gamma}_{k}, e_{k-1}}\|f_{\gamma_{k}, e_{k-1}})}{\partial
\rho}<0$.

\begin{figure}[htbp]
\centering
\includegraphics[height=0.8\columnwidth,width=\columnwidth]{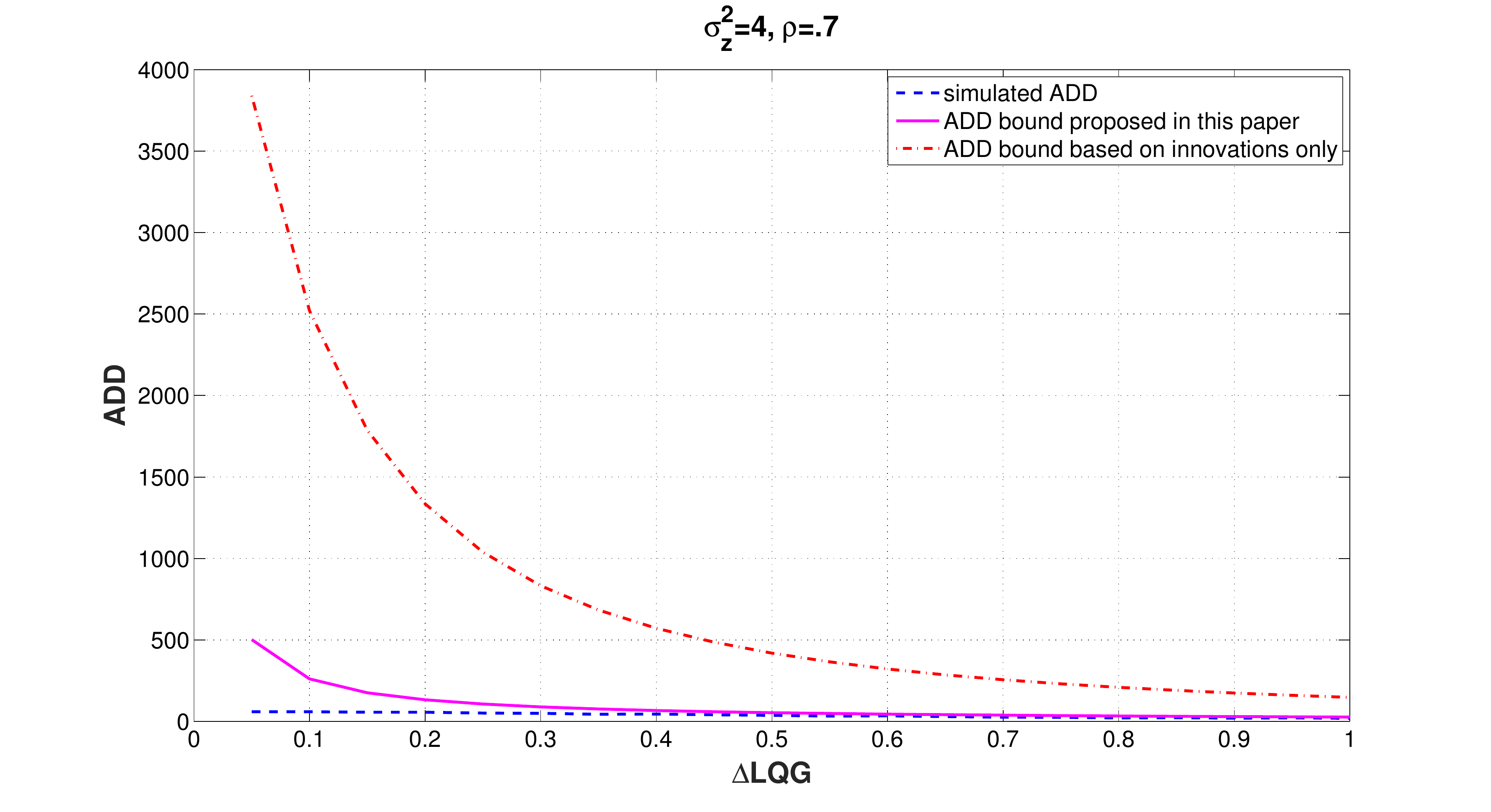}
\\
\includegraphics[height=0.8\columnwidth,width=\columnwidth]{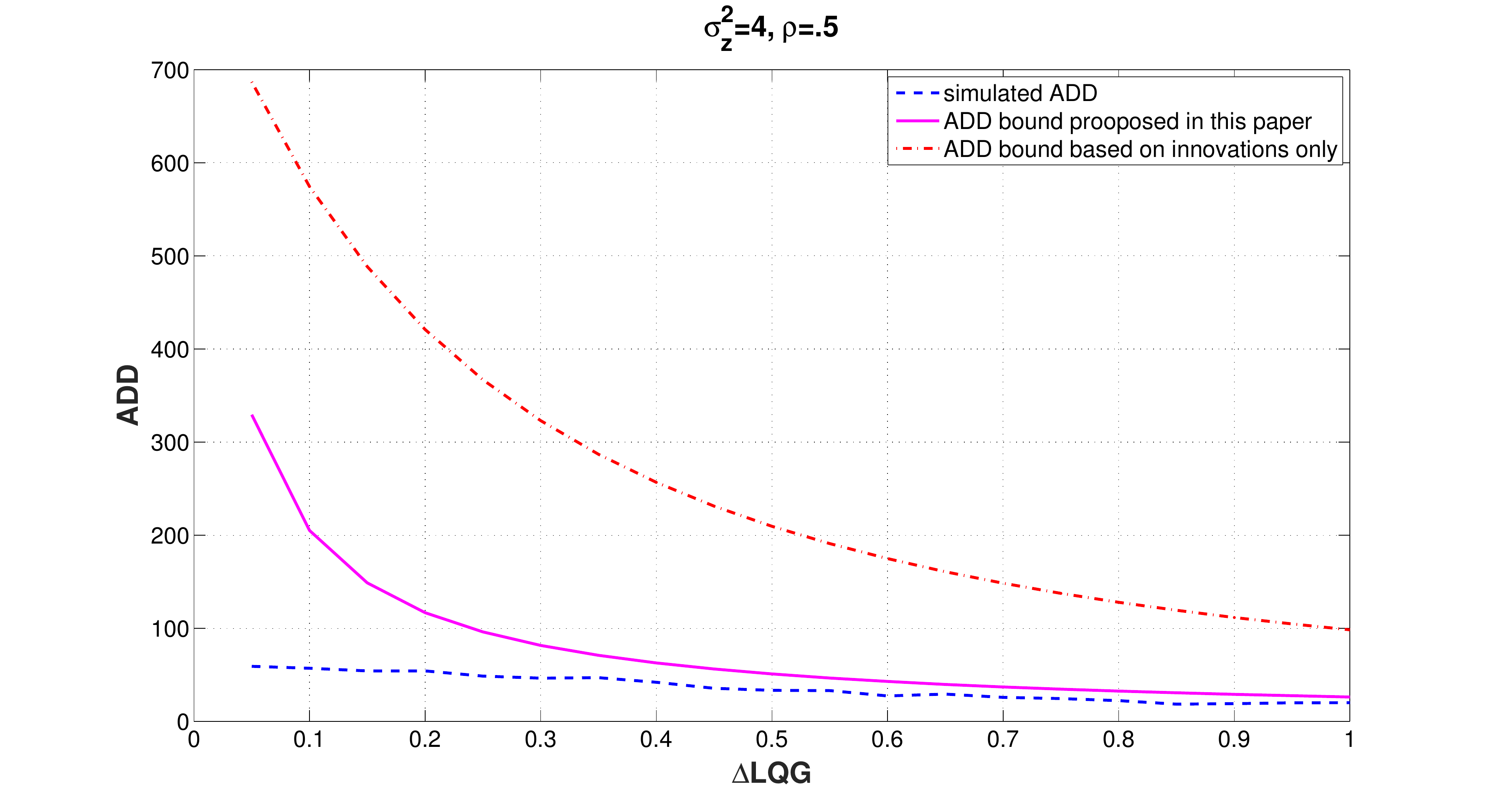}
 \caption{\footnotesize{\footnotesize{ADD in terms of $\Delta LQG$ for different values of $\sigma_{z}^{2}$ and $\rho$ with $A=1.2$}} }
\label{fig3}
\end{figure}

Fig.~\ref{fig3} illustrates the ADD for an unstable system with $A>1$. The ADD is shown in two different
setup for $A=1.2$ while $B,C,R,Q,W,U,p^{F}$ are the same as in
Fig.~\ref{fig2}. Comparing the corresponding diagrams, i.e., the same
$\sigma_{z}^{2}$ and $\rho$ in Fig. \ref{fig2} and Fig. \ref{fig3}, it is seen  that in
the unstable case with $A>1$, there exists a higher gap between the ADD bounds for the sequential tests based on the joint distribution as proposed in this paper and
the one based on innovations/residues only.

\section{Conclusions and Future Work}
In this paper, we investigated how a suitably designed sequential detection test can detect deception attacks in a scalar networked control system with an average detection delay that can be reduced by introducing a physical watermarking signal with a suitable variance. The tradeoff between quick detection and penalty in the control cost as a result of using the watermarking signal is also investigated. Future works will extend these results to multi-variable (vector state and measurements) systems, and propose a  dynamic game between the adversary and the control system designer regarding the attacker's effort to remain stealthy, and the system designer's effort to detect the attack with minimum delay, using such physical watermarking schemes. 
\appendix
\subsection{Proof of Theorem 1}
To calculate $D(f_{\widetilde{\gamma}_{k}, e_{k-1}}\|f_{\gamma_{k}, e_{k-1}})$, we need to obtain the joint
distributions $f_{\gamma_{k}, e_{k-1}}(\gamma,e)$ and
$f_{\widetilde{\gamma}_{k}, e_{k-1}}(\gamma,e)$ in steady state. As it
was shown in \eqref{eq1}, in the healthy system, $\gamma_{k}$ and
$e_{k-1}$ are uncorrelated for i.i.d. watermarking sequence. Hence
the joint distribution $f_{\gamma_{k}e_{k-1}}(\gamma,e)$ appears as:

\begin{eqnarray}
f_{\gamma_{k}e_{k-1}}(\gamma,e)=\frac{1}{2\pi
\sigma_{\gamma}\sigma_{e}}\exp\frac{-1}{2}(\frac{\gamma^{2}}{\sigma_{\gamma}^{2}}+\frac{e^{2}}{\sigma_{e}^{2}})\label{noattackdis}
\end{eqnarray}
in which $\sigma_{\gamma}^{2}$ is given as in
\eqref{sigmagammanoattack}. On the other hand, when we have attack,
$\tilde{\gamma}_{k}$ and $e_{k-1}$ are correlated according to
\eqref{eq2}. Since $\tilde{\gamma}_{k}$ and $e_{k-1}$ are zero-mean
Gaussian, to obtain their joint distribution, we need to calculate
$\sigma_{\tilde{\gamma}}^{2}$ and
$\textbf{cov}(\tilde{\gamma}_{k},e_{k-1})$. Since $e_{k-1}$ is
uncorrelated with $z_{k}$ and $\hat{x}^F_{k-1|k-1}$, we use \eqref{eq2} to
obtain:
\begin{equation}
\textbf{cov}(\tilde{\gamma}_{k},e_{k-1})=-CB\sigma_{e}^{2}\label{eq7}
\end{equation}
Using the same equation, we have:
\begin{align}
\sigma_{\tilde{\gamma}}^{2} &=\sigma_{z}^{2}+C^{2}(A+BL)^{2}\sigma_{\hat{x}}^{2} 
\nonumber \\
& -2C(A+BL)\textbf{cov}(z_{k},\hat{x}^F_{k-1|k-1})
+C^{2}B^{2}\sigma_{e}^{2}.\label{eq6}
\end{align}
where $\sigma_{\hat{x}}^{2} = \mathrm{E}\left(\hat{x}^F_{k-1|k-1}\right)^2$. 

To calculate $\textbf{cov}(z_{k},\hat{x}^F_{k-1|k-1})$ we proceed as
follows. By combining \eqref{xestimatedef} and \eqref{gammadef}, for
the attacked system, one obtains:
\begin{equation}
\hat{x}^F_{k-1|k-1}=Kz_{k-1}+\mathcal{A}\hat{x}^F_{k-2|k-2}+B(1-CK)e_{k-2}.\label{eq3}
\end{equation}

By multiplying the above equation with $z_{k}$ and calculating
expectation of the both sides for the stationary system, and defining $\mathrm{E}_{\hat{x}z}(-l) = 
\textbf{cov}(z_{k},\hat{x}^F_{k-l|k-l})$,
we have:
\begin{equation}
\mathrm{E}_{\hat{x}z}(-1)=KE_{zz}(1)+\mathcal{A}E_{\hat{x}z}(-2).\no
\end{equation}
Continuing the same procedure for $z_{k+1},z_{k+2},...$, one
obtains:
\begin{eqnarray}
&\mathrm{E}_{\hat{x}z}(-2)=K\mathrm{E}_{zz}(2)+\mathcal{A}\mathrm{E}_{\hat{x}z}(-3)\no\\
&\mathrm{E}_{\hat{x}z}(-3)=K\mathrm{E}_{zz}(3)+\mathcal{A}\mathrm{E}_{\hat{x}z}(-4)\no\\
&.\no\\
&.\no\\
&.\no
\end{eqnarray}
in which $\mathrm{E}_{zz}(k)$ is obtained according to
\eqref{zcorrelation}.

Since $\mathcal{A}<1$ and $\rho<1$, $\mathrm{E}_{\hat{x}z}(-1)$ is
obtained as a sum of an infinite geometric series which converges
to:
\begin{eqnarray}
&\mathrm{E}_{\hat{x}z}(-1)=K\sigma_{z}^{2}\sum\limits_{i=1}^{\infty}\rho^{i}\mathcal{A}^{i-1}\no\\
&=K\sigma_{z}^{2}\rho\sum\limits_{i=0}^{\infty}(\rho\mathcal{A})^{i}\no\\
&=\frac{K\sigma_{z}^{2}\rho}{1-\rho\mathcal{A}}\label{eq4}
\end{eqnarray}

To calculate $\sigma_{\hat{x}}^{2}$, we reuse \eqref{eq3}
such that:
\begin{align}
\mathrm{E}(\hat{x}^F_{k-1|k-1})^2 & =K^{2}\mathrm{E}(z_{k-1}^{2})+\mathcal{A}^{2}\mathrm{E}(\hat{x}^F_{k-2|k-2})^{2}
\nonumber \\
& +2K\mathcal{A}\mathrm{E}(z_{k-1}\hat{x}^F_{k-2|k-2})
+B^{2}(1-CK)^{2}\sigma_{e}^{2}
\end{align}
which results in:
\begin{equation}
\sigma_{\hat{x}}^2=K^{2}\sigma_{z}^2+\mathcal{A}^{2}\sigma_{\hat{x}}^2+2K\mathcal{A}\mathrm{E}_{\hat{x}z}(-1)+B^{2}(1-CK)^{2}\sigma_{e}^{2}.\label{eq5}
\end{equation}
Combing \eqref{eq4} and \eqref{eq5} yields:
\begin{equation}
\sigma_{\hat{x}}^2=\frac{K^{2}(1+\rho\mathcal{A})}{(1-\rho\mathcal{A})(1-\mathcal{A}^{2})}\sigma_{z}^2+\frac{B^{2}(1-CK)^{2}}{1-\mathcal{A}^{2}}\sigma_{e}^{2}.
\end{equation}

and eventually, using \eqref{eq6}, $\sigma_{\tilde{\gamma}}^{2}$ is
obtained as:
\begin{eqnarray}
&\sigma_{\tilde{\gamma}}^{2}=[(1-\frac{\rho
CK(A+BL)}{1-\rho\mathcal{A}})^{2}+\frac{(1-\rho^{2})C^{2}K^{2}(A+BL)^{2}}{(1-\mathcal{A}^{2})(1-\rho\mathcal{A})^{2}}]\sigma_{z}^{2}\no\\
&+\frac{B^{2}C^{2}}{1-\mathcal{A}^{2}}\sigma_{e}^{2}\no
\end{eqnarray}

To obtain the joint distribution
$f_{\widetilde{\gamma}_{k}e_{k-1}}(\gamma,e)$, we form the
cross-covariance matrix of $\gamma_{k}.e_{k-1}$ as:
\begin{eqnarray}
\Sigma=
\left(%
\begin{array}{cc}
  \sigma_{\tilde{\gamma}}^{2} & -CB\sigma_{e}^{2} \\
  -CB\sigma_{e}^{2} & \sigma_{e}^{2} \\
\end{array}%
\right)
\end{eqnarray}
and consequently:
\begin{align}
f_{\widetilde{\gamma}_{k}, e_{k-1}}(\gamma,e) & = \frac{1}{2\pi
\sigma_{\tilde{\gamma}}\sigma_{e}\sqrt{1-\lambda^{2}}}
\nonumber \\
& \times \exp \left\{ \frac{-1}{2(1-\lambda^{2})}\left(\frac{\gamma^{2}}{\sigma_{\tilde{\gamma}}^{2}}+\frac{e^{2}}{\sigma_{e}^{2}}-\frac{2\lambda
e\gamma}{\sigma_{e}\sigma_{\tilde{\gamma}}}\right) \right\} \label{attackdis}
\end{align}
in which
\begin{equation}
\lambda=\frac{\textbf{cov}(\tilde{\gamma}_{k},e_{k-1})}{\sigma_{e}\sigma_{\tilde{\gamma}}}\mathop{=}\limits^{(a)}\frac{-BC\sigma_{e}}{\sigma_{\tilde{\gamma}}}
\end{equation}
where (a) is deduced from \eqref{eq7}. Note also that $|\lambda| < 1$ as it is a correlation coefficient.

Then $D(f_{\widetilde{\gamma}_{k}, e_{k-1}} \|f_{\gamma_{k}, e_{k-1}})$ is calculated as:
\begin{eqnarray}
\iint_{-\infty}^{\infty}f_{\widetilde{\gamma}_{k}e_{k-1}}(\gamma,e)\log\frac{f_{\widetilde{\gamma}_{k},e_{k-1}}(\gamma,e)}{f_{\gamma_{k},e_{k-1}}(\gamma,e)})de
d\gamma.\label{eq8}
\end{eqnarray}

Replacing the joint distributions as in \eqref{noattackdis} and
\eqref{attackdis} in \eqref{eq8} yields
\begin{align}
&D(f_{\widetilde{\gamma}_{k}, e_{k-1}} \|f_{\gamma_{k}, e_{k-1}}) = 
\log(\frac{\sigma_{\gamma}}{\sigma_{\tilde{\gamma}}\sqrt{1-\lambda^{2}}}) \no \\
& +\frac{-\sigma_{\tilde{\gamma}}^{2}}{2}
(\frac{1}{(1-\lambda^{2})\sigma_{\tilde{\gamma}}^{2}}-\frac{1}{\sigma_{\gamma}^{2}})+\frac{-1}{2}
(\frac{1}{(1-\lambda^{2})}-1)\no\\
& \;\; \; +\frac{\lambda
\textbf{cov}(\tilde{\gamma}_{k},e_{k-1})}{(1-\lambda^{2})\sigma_{\tilde{\gamma}}\sigma_{e}}\no\\
&=\log(\frac{\sigma_{\gamma}}{\sigma_{\tilde{\gamma}}\sqrt{1-\lambda^{2}}})+\frac{-\sigma_{\tilde{\gamma}}^{2}}{2}
(\frac{1}{(1-\lambda^{2})\sigma_{\tilde{\gamma}}^{2}}-\frac{1}{\sigma_{\gamma}^{2}}) \no \\ 
& +\frac{-1}{2}
(\frac{1}{(1-\lambda^{2})}-1)
+\frac{\lambda^{2}}{(1-\lambda^{2})}\no\\
&=\log(\frac{\sigma_{\gamma}}{\sigma_{\tilde{\gamma}}\sqrt{1-\lambda^{2}}})+\frac{\sigma_{\tilde{\gamma}}^{2}}{2\sigma_{\gamma}^{2}}-\frac{1}{2}\no\\
&=\frac{1}{2}\log(\frac{1}{1-\lambda^{2}})+
\frac{1}{2}(\frac{\sigma_{\tilde{\gamma}}^{2}}{\sigma_{\gamma}^{2}}-1-\log\frac{\sigma_{\tilde{\gamma}}^{2}}{\sigma_{\gamma}^{2}})\no
\end{align}
Using the same approach as in \cite{enoch}, it can be shown that the
averaged KLD is equal to the single letter $D(f_{\widetilde{\gamma}_{k}, e_{k-1}} \|f_{\gamma_{k}, e_{k-1}})$ calculated
above in which the distributions are the steady state ones.

\subsection{Proof of Theorem 2}
To calculate the cost of using the watermarking for the purpose of
detection, we obtain the difference of LQGs between the cases depending on whether  watermarking
is used or not used  in the healthy system. According to \eqref{LQG}, we use \eqref{uwithw} to
obtain:
\begin{eqnarray}
&\Delta LQG=J_{w}-J_{n}=\no\\
& W(\mathrm{E}(X_{w}^{2})\!-\!\mathrm{E}(X_{n}^{2}))\!+\!
UL^{2}(\mathrm{E}(\hat{X}_{w}^{2})\!-\!\mathrm{E}(\hat{X}_{n}^{2}))\!+\!U\sigma_{e}^{2}\label{deltaLQG}
\end{eqnarray}
where subscript 'w' refers to the case where we use watermarking and
subscript 'n' refers to the case that we don't use watermarking. To
calculate $\Delta LQG$, we need to calculate $\mathrm{E}(X_{w}^{2})$
and $\mathrm{E}(\hat{X}_{w}^{2})$. Based on the fact that
$\gamma_{k}$ is uncorrelated with $\hat{x}_{k|k-1}$, we use
\eqref{xestimatedef} in steady-state condition to obtain:

\begin{eqnarray}
&\mathrm{E}(\hat{X}_{w}^{2})=(A+BL)^{2}\mathrm{E}(\hat{X}_{w}^{2})+B^{2}\sigma_{e}^{2}+K^{2}\sigma_{\gamma}^{2}\no\\
&=(A+BL)^{2}\mathrm{E}(\hat{X}_{w}^{2})+B^{2}\sigma_{e}^{2}+K^{2}(C^{2}P+R)\no
\end{eqnarray}
which yields
\begin{equation}
\mathrm{E}(\hat{X}_{w}^{2})=\frac{1}{1-(A+BL)^{2}}(B^{2}\sigma_{e}^{2}+K^{2}(C^{2}P+R))\label{expxhat2}
\end{equation}
To calculate $\mathrm{E}(X_{w}^{2})$, \eqref{gammadef} is obtained
as:
\begin{eqnarray}
&\mathrm{E}(Y_{w}^{2})=\sigma_{\gamma}^{2}+C^{2}((A+BL)^{2}\mathrm{E}(\hat{X}_{w}^{2})+B^{2}\sigma_{e}^{2})\no\\
&=(C^{2}P+R)(1+\frac{C^{2}K^{2}(A+BL)^{2}}{1-(A+BL)^{2}})+\frac{B^{2}C^{2}}{1-(A+BL)^{2}}\sigma_{e}^{2})\label{eq10}
\end{eqnarray}
where we used \eqref{expxhat2} to conclude \eqref{eq10}.

Then \eqref{ydef} is used to calculate $\mathrm{E}(X_{w}^{2})$ as
\begin{eqnarray}
&\mathrm{E}(X_{w}^{2})=\frac{\mathrm{E}(Y_{w}^{2})-R}{C^{2}}\no\\
&=P+\frac{K^{2}(A+BL)^{2}(C^{2}P+R)}{1-(A+BL)^{2}}+\frac{B^{2}}{1-(A+BL)^{2}}\sigma_{e}^{2}\label{eq11}
\end{eqnarray}

Combining \eqref{expxhat2} and \eqref{eq11} with \eqref{deltaLQG}
gives us:
\begin{eqnarray}
\Delta
LQG=\bigg(U+\frac{B^{2}(W+L^{2}U)}{1-(A+BL)^{2}}\bigg)\sigma_{e}^{2}
\end{eqnarray}

\end{document}